\documentclass[12pt,a4paper,twoside,final,notitlepage, leqno]{article}
\usepackage[english]{babel}
\usepackage[T1]{fontenc}
\usepackage{epsfig, graphicx, amssymb}
\usepackage{amsmath,amsthm,epsfig,amsfonts}
\usepackage{float}
\usepackage{color}
\def\br {\break}
\newcommand{\moneq}{\vspace*{-7pt} \begin{equation} \displaystyle }
\newcommand{\moneqstar}{\vspace*{-6pt} \begin{equation*} \displaystyle }
\newcommand{\monendstar}{\vspace*{-6pt} \end{equation*}   }
\newcommand{\monend}{\vspace*{-7pt} \end{equation}   }
\newcommand{\moneqarraystar}{ \begin{eqnarray*} \displaystyle }
\newcommand{\monendarraystar}{ \end{eqnarray*}   }
\newcommand{\dd}{{\rm d}}

\definecolor{vertfonce}{rgb}{0.0, 0.5, 0.0}
\def\section*#1{}
\usepackage{fancyhdr}
\renewcommand{\headrulewidth}{0pt}
\begin{document}

\fancypagestyle{plain}{ \fancyfoot{} \renewcommand{\footrulewidth}{0pt}}
\fancypagestyle{plain}{ \fancyhead{} \renewcommand{\headrulewidth}{0pt}}
%%   \bibliographystyle{alpha}

%%%%%%%%%%%%%%%%%%%%%%%%%%%%%%%%%%%%%%%%%%%%%%%%%%%%%%%%%%%%%%%%%%%%%%%%%%%%%%%
~

  \vskip 2.1 cm

\centerline {\bf \LARGE A variational  symplectic scheme}

\bigskip

\centerline {\bf \LARGE based on Simpson's quadrature}

 \bigskip  \bigskip \bigskip

\centerline { \large    Fran\c{c}ois Dubois$^{ab}$ and Juan Antonio Rojas-Quintero$^{c}$}

\smallskip  \bigskip

\centerline { \it  \small
  $^a$   Laboratoire de Math\'ematiques d'Orsay, Facult\'e des Sciences d'Orsay,}

\centerline { \it  \small   Universit\'e Paris-Saclay, France.}

\centerline { \it  \small
$^b$    Conservatoire National des Arts et M\'etiers, LMSSC laboratory,  Paris, France.}

\centerline { \it  \small
$^c$    CONAHCYT/Tecnol\'ogico Nacional de M\'exico/I.T. Ensenada, Ensenada 22780, BC, Mexico.}

%%%   \bigskip

\bigskip  \bigskip

\centerline {11 May 2023 
{\footnote {\rm  \small $\,$ This contribution has been presented to the
6th conference on Geometric Science of Information,
%% ``From classical to quantum information geometry'',
Saint-Malo (France), 30 august - 01 september 2023. Published in
{\it Geometric Science of Information},  LNCS 14072, 
%%  6th International Conference, GSI 2023, St. Malo, France, August 30 – September 1, 20
F. Nielsen and F. Barbaresco (Eds.), pages 22-31, 2023.
Edition 22 June 2024.}}}
%%%    {\footnote



%%%  76N15,  %%%  Gas dynamics, general
%%%  82C20.   %%%  Dynamic lattice systems (kinetic Ising, etc.) and systems on graphs

\bigskip  \bigskip
\noindent {\bf \large Abstract}

\noindent
We propose a variational  symplectic numerical method for the time integration of
dynamical systems issued from the least action principle.
We  assume a quadratic internal interpolation of the state and we approximate the action
in a small time step by the Simpson's quadrature formula.
The resulting scheme is explicited for an elementary harmonic oscillator. It is a stable, 
explicit, and  symplectic scheme satisfying the conservation of an approximate energy.
Numerical tests illustrate our theoretical study. 

%%%%%%%%%%%%%%%%%%%%%%%%%%%%%%%%%%%%%%%%%%%%%%%%%%%%%%%%%%%%%%%%%%%%%%%%%%%%%%%  section 1
\bigskip \bigskip    \noindent {\bf \large    1) \quad  Introduction} 
%%%%%%%%%%%%%%%%%%%%%%%%%%%%%%%%%%%%%%%%%%%%%%%%%%%%%%%%%%%%%%%%%%%%%%%%%%%%%%%%%%%%%%%%%%

\smallskip \noindent
The principle of least action is a key point for establishing evolution equations
or partial differential equations, from classical to quantum mechanics and electromagnetisms
\cite{Ar74,FH65,So70}. An important application of this principle is proposed
with the finite element method \cite{Co43} and it is used for  engineering applications 
since the 1950's.
For dynamics equations and dynamical systems, 
a synthesis of the state of the art is proposed in \cite{HLW06,SS92,WM97}.

\smallskip \noindent
In this contribution, we first recall the classical variational approach.
It is founded on  a midpoint %%% trapezoidal
quadrature formula for the approximate calculation of an integral. 
We essentially follow the contribution \cite{KMOW00}
in this Section 2. 
Then we recall in Section 3 the interpolation of functions with quadratic finite elements.
Once this prerequisite is  in place, we develop in Section 4 the approximation
of discrete Lagrangians with Simpson's quadrature formula.
The result is a numerical scheme that can be considered as a variant of the classical
approach presented in Section 2 and we derive in Section 5 the discrete Euler-Lagrange equations.
We notice in Section 6 that the scheme admits a symplectic structure and in Section 7
that an approximation of the
energy is conserved along the discrete time integration.
First numerical results are presented in Section 8 before some words of conclusion.

%%%%%%%%%%%%%%%%%%%%%%%%%%%%%%%%%%%%%%%%%%%%%%%%%%%%%%%%%%%%%%%%%%%%%%%%%%%%%%%  section 2
\bigskip    \noindent {\bf \large    2) \quad  A classical  variational symplectic numerical scheme} 
%%%%%%%%%%%%%%%%%%%%%%%%%%%%%%%%%%%%%%%%%%%%%%%%%%%%%%%%%%%%%%%%%%%%%%%%%%%%%%%%%%%%%%%%%%

\smallskip \noindent
We consider a dynamical system described by a state $ \, q(t) \, $ composed by a simple real variable 
to fix the ideas,  and for $ \, 0 \leq t \leq T $.
The continuous action~$ \, S_c \, $ introduces a Lagrangian~$ \, L \, $
and we have
\moneq \label{action-continue} 
S_c = \int_0^T L \Big( {{\dd q}\over{\dd t}} ,\, q(t) \Big) \, \dd t . 
\monend  
We use in this contribution a very classical Lagrangian 
\moneq \label{Lagrangien-continu} 
L \Big( {{\dd q}\over{\dd t}} ,\, q \Big)  = {m\over2} \, \Big( {{\dd q}\over{\dd t}} \Big)^2 - V(q) . 
\monend  
A discretization of the relation (\ref{action-continue}) is obtained  by splitting the interval
$ \, [0 ,\, T ] \, $  into $ \, N \, $ elements and we set $\, h = {{T}\over{N}} $.
At the discrete time $ \, t_j = j \, h $, an approximation  $\, q_j \, $ of
$ \, q(t_j) \, $ is introduced and a discrete form $ \, S_d \, $  of the continuous action $ \, S_c \, $ 
can be defined according to
\moneqstar 
 S_d = \sum_{j=1}^{N-1} L_d (q_j ,\, q_{j+1} ) .
\monendstar
The discrete Lagrangian  $ \, L_d (q_\ell ,\, q_r ) \, $
is derived from the  relation  (\ref{Lagrangien-continu})  
with a centered  finite difference approximation 
$ \,  {{\dd q}\over{\dd t}} \simeq {{q_r - q_\ell}\over{h}} $ 
and a midpoint  quadrature formula
%
%
%%%%%%%%%%%%%%%%%%%%%%%%%%%%%%%%%%%%%%%%%%%%%%%%%%%%%%%%%%%%%%%%%%%%%%%%%%%%%%%%%%%%%%%%%
\fancyhead[EC]{\sc{Fran\c{c}ois Dubois and Juan Antonio Rojas-Quintero}}
\fancyhead[OC]{\sc{A variational  symplectic scheme based on Simpson's quadrature}}
%%%%%%%%%%%%%%%%%%%%%%%%%%%%%%%%%%%%%%%%%%%%%%%%%%%%%%%%%%%%%%%%%%%%%%%%%%%%%%%%%%%%%%%%%
%
%
\moneqstar 
\int_0^h V \big( q(t) \big) \, \dd t \simeq h \, V \Big( {{q_\ell + q_r}\over2} \Big)  \, :
\monendstar
%
%%%%%%%%%%%%%%%%%%%%%%%%%%%%%%%%%%%%%%%%%%%%%  jolie numerotation des pages
\fancyfoot[C]{\oldstylenums{\thepage}}
%%%%%%%%%%%%%%%%%%%%%%%%%%%%%%%%%%%%%%%%%%%%%  fin jolie numerotation des pages
%
\vskip -.3cm 
\moneq \label{lagrangien-discret} 
%%% L_d (q_\ell ,\, q_r ) = {{m\,h}\over2} \, \Big( {{q_r - q_\ell}\over{h}} \Big)^2 - {{h}\over2} \, \big( V (q_r) +  V (q_\ell)  \big)  .
L_d (q_\ell ,\, q_r ) = {{m\,h}\over2} \, \Big( {{q_r - q_\ell}\over{h}} \Big)^2 - h \, V \Big( {{q_\ell + q_r}\over2} \Big) .
\monend
\noindent 
We observe that  $ \, S_d = \cdots  + \,  L_d (q_{j-1} ,\, q_j ) +  L_d (q_j ,\, q_{j+1} )  \, + \cdots $.
Then the discrete Euler Lagrange equation $ \, \delta S_d = 0 \, $ for  an arbitrary variation $ \, \delta q_j \, $
of the discrete variable $ \, q_j \, $ can be written
\moneq \label{Euler-Lagrange-discret} 
{{\partial L_d}\over{\partial q_r}}  (q_{j-1} ,\, q_j )  + {{\partial L_d}\over{\partial q_\ell}}  (q_j ,\, q_{j+1} ) = 0 .
\monend
Taking into account the  relation (\ref{lagrangien-discret}), we obtain
\moneq   \label{schema-ordre-deux} 
m \, {{q_{j+1} - 2\, q_j + q_{j-1}}\over{h^2}} + {1\over2} \, \Big[ {{\dd V}\over{\dd q}} \Big( {{q_j + q_{j+1}}\over2} \Big)
+  {{\dd V}\over{\dd q}} \Big( {{q_{j-1} + q_{j}}\over2} \Big)  \Big] = 0 . 
\monend 
This numerical scheme is clearly consistent with the second order differential
equation
\moneq   \label{oscillateur-nonlineaire} 
m \, {{\dd^2 q}\over{\dd t^2}} + {{\dd V}\over{\dd q}} = 0 
\monend 
associated with the Lagrangian
proposed in (\ref{Lagrangien-continu}). It is easy to verify that when 
\moneq \label{oscillateur harmonique} 
V(q) = {1\over2} \, m\, \omega^2 \, q^2  ,
\monend
the scheme (\ref{schema-ordre-deux}) is linearly stable.
We suppose that the assumption  (\ref{oscillateur harmonique})
is satisfied until the end of this paragraph. 
The momentum $ \, p_r \, $ is defined by
\moneq \label{moment} 
p_r = {{\partial L_d}\over{\partial q_r}}  (q_\ell  ,\, q_r ) .
\monend
We have $ \, p_{j+1} = m \, {{q_{j+1} - q_j}\over{h}} - m \,{{\omega^2 \, h}\over{4}} \, \big( q_{j+1} + q_j \big) \, $
and an analogous relation for $ \, p_j $. Then after some lines of algebra, we obtain
a discrete system involving the momentum and the state:
\moneq \label{systeme-dynamique-discret-ordre-2} 
p_{j+1} - p_j = - m \, {{\omega^2 \, h}\over{2}} \big( q_{j+1} + q_j \big) \,,\,\,
{{q_{j+1} - q_j}\over{h}} = {1\over{2\, m}} \,  \big( p_{j+1} + p_j \big) .
\monend
These relations are consistent with the first order Hamilton version
$ \, {{\dd p}\over{\dd t}} +  m \,\omega^2 \, q = 0 $, $ \,  {{\dd q}\over{\dd t}} = {{p}\over{m}} \, $ 
of the equations of an harmonic oscillator.
Moreover, we can write the system (\ref{systeme-dynamique-discret-ordre-2})
under the form
\moneq \label{systeme-symplectique-ordre-2}
\begin{pmatrix}p_{j+1} \\ q_{j+1} \end{pmatrix} =  \Phi \,  \begin{pmatrix}p_{j} \\ q_{j} \end{pmatrix} 
\monend
with
\moneq \label{matrice-symplectique-ordre-2}
\Phi = {{1}\over{1+{{\omega^2 \, h^2}\over4}}}
\begin{pmatrix}  1-{{\omega^2 \, h^2}\over4} & - m \, \omega^2 \, h \\ {{h}\over{m}} &  1-{{\omega^2 \, h^2}\over4} \end{pmatrix} .
\monend
Because $\,  \det \Phi  = 1 $, the discrete flow (\ref{systeme-symplectique-ordre-2}) is symplectic
as observed by Sanz-Serna \cite{SS92}.
Moreover,  Kane {\it et al.} \cite{KMOW00} have  remarked that 
the numerical scheme (\ref{systeme-symplectique-ordre-2}) 
is one particular inconditionally stable version of the
Newmark scheme \cite{Ne59}. Last but not least,
the discrete Hamiltonian~$\, H_j \, $ defined by 
\moneq \label{HH-exact}  
H_j \equiv {1\over{2 \,m}}\, p_j^2 + {{m \, \omega^2}\over{2}} \, q_j^2
\monend
is conserved: we have $ \, H_{j+1} = H_j \,$ for $ \, 0 \leq j \leq N-1 $. 
%%  \noindent
We consider now a more elaborate interpolation in each interval,
updating affine functions by polynomials of degree two.

%%%%%%%%%%%%%%%%%%%%%%%%%%%%%%%%%%%%%%%%%%%%%%%%%%%%%%%%%%%%%%%%%%%%%%%%%%%%%%%  section 3
\bigskip    \noindent {\bf \large    3) \quad  Quadratic interpolation} %%  finite element} 
%%%%%%%%%%%%%%%%%%%%%%%%%%%%%%%%%%%%%%%%%%%%%%%%%%%%%%%%%%%%%%%%%%%%%%%%%%%%%%%%%%%%%%%%%%

\smallskip \noindent
Internal interpolation between $ \, 0 \, $ and $ \, h \, $ can be written
in terms of quadratic finite elements~\cite{RT83}. For $ \, 0 \leq \theta \leq 1 $, we first set 
\moneq \label{fonctions-de-base} 
\varphi_0(\theta) =  (1 - \theta) \, (1 - 2 \, \theta) \,,\,\,  
\varphi_{1/2}(\theta) =  4 \, \theta  \, (1 - \theta) \,,\,\,   
\varphi_{1}(\theta) =  \theta \, (2 \, \theta  - 1 ) .   
\monend
With $ \, t = h \, \theta $, we consider the polynomial function
\moneq \label{q2t} 
q(t) = q_\ell \,\varphi_0(\theta) + q_m \, \varphi_{1/2}(\theta) +  q_r \, \varphi_{1}(\theta) . 
\monend 
Then $ \, q(0) = q_\ell $, $ \, q({{h}\over2}) =  q_m \, $ and $ \, q(h) = q_r \, $
and the basis functions (\ref{fonctions-de-base}) are well adapted to these degrees of freedom. 
We have also
\moneqstar    \begin{array} {l}
{{\dd q}\over{\dd t}} = {{1}\over{h}} \, \big[ q_\ell \, {{\dd \varphi_0}\over{\dd \theta}}
+  q_m \, {{\dd \varphi_{1/2}}\over{\dd \theta}} +  q_r \,   {{\dd \varphi_1}\over{\dd \theta}} \big]
\\   \quad \,  =
      {{1}\over{h}} \, \big[  q_\ell \, (4 \, \theta - 3) + 4 \,  q_m \, (1 - 2 \, \theta) + q_r \, ( 4 \, \theta - 1) \big]
\\   \quad \,  = g_\ell \, (1-\theta) + g_r \, \theta 
\end{array}  \monendstar

%%%     \noindent 
%%%     $  {{\dd q}\over{\dd t}} = {{1}\over{h}} \, \big[ q_\ell \, {{\dd \varphi_0}\over{\dd \theta}}
%%%     +  q_m \, {{\dd \varphi_{1/2}}\over{\dd \theta}} +  q_r \,   {{\dd \varphi_1}\over{\dd \theta}} \big] $
%%%     $ \,\, =
%%%     {{1}\over{h}} \, \big[  q_\ell \, (4 \, \theta - 3) + 4 \,  q_m \, (1 - 2 \, \theta) + q_r \, ( 4 \, \theta - 1) \big] $ 
%%%     \smallskip \noindent \quad $ \,\, = g_\ell \, (1-\theta) + g_r \, \theta \, $

\noindent
 with the derivatives $ \, g_\ell \, $ and $ \, g_r \, $ given by a Gear scheme \cite{Ge71}, {\it id est} 
\moneq \label{derivees-gear} 
g_\ell =  {{\dd q}\over{\dd t}}(0) =  {{1}\over{h}} \,  \big( -3 \, q_\ell + 4 \,  q_m - q_r \big) \,,\,\, 
g_r =  {{\dd q}\over{\dd t}}(h) = {{1}\over{h}} \,  \big(  q_\ell - 4 \,  q_m + 3 \, q_r \big) .
\monend 
We remark also that 
\moneq \label{derivee-centree}
 g_m = {{\dd q}\over{\dd t}} \Big( {{h}\over2} \Big) = {1\over2} \, (g_\ell + g_r) = {{q_r - q_\ell}\over{h}} .
\monend 

\smallskip \noindent 
Once the interpolation is defined in an interval of length $ \, h $, we use it by splitting
the range~$ \, [0,\,T] \, $  into $ \, N \, $  pieces, and $ \, h = {{T}\over{N}} $.
With $ \, t_j = j \, h $, we set $ \, q_j \simeq q(t_j) \, $ for $ \, 0 \leq j \leq N \, $
and $ \, q_{j+1/2} \simeq q(t_j+{{h}\over{2}}) \, $ with $ \, 0 \leq j \leq N-1 $.
In the interval $ \, [t_j ,\, t_{j+1} ] $, the function $ \, q(t) \, $ is a polynomial of degree 2, represented by the
relation (\ref{q2t}) with $ \, t = t_j + \theta \, h  $, $ \, q_\ell = q_j $, $\, q_m = q_{j+1/2} \, $
and $ \, q_{r} =  q_{j+1} $.

%%%%%%%%%%%%%%%%%%%%%%%%%%%%%%%%%%%%%%%%%%%%%%%%%%%%%%%%%%%%%%%%%%%%%%%%%%%%%%%  section 4
\bigskip    \noindent {\bf \large    4) \quad    Simpson's quadrature  for a discrete Lagrangian} 
%%%%%%%%%%%%%%%%%%%%%%%%%%%%%%%%%%%%%%%%%%%%%%%%%%%%%%%%%%%%%%%%%%%%%%%%%%%%%%%%%%%%%%%%%%

\smallskip \noindent
For the numerical integration of a regular function $ \, \psi \, $ on the interval $\, [0,\, 1] $,
the midpoint  %%% trapezoidal
method studied previously
%%% $ \, \int_0^1 \psi(\theta) \, \dd \theta \simeq {1\over2} \, \big( \psi(0) + \psi(1) \big) \, $ 
$ \, \int_0^1 \psi(\theta) \, \dd \theta \simeq \psi \big( {1\over2} \big) \, $ 
is exact for a polynomial $ \, \psi \, $ of degree smaller or equal to $ \, 1 $.  %%%  $ \, 2 $. 
To obtain a better precision, a very popular method has been proposed by Thomas Simpson (1710-1761):
\moneq \label{simpson} 
 \int_0^1 \psi(\theta) \, \dd \theta \simeq {1\over6} \, \Big[ \psi(0) + 4 \, \psi\Big( {1\over2} \Big) + \psi(1) \Big] .
\monend 
The quadrature formula  (\ref{simpson})  is accurate up to polynomials of  degree three. 
Then a discrete  Lagrangian
$ \, L_h (q_\ell ,\, q_m ,\, q_r) \simeq \int_0^h  \big[ {m\over2}  \big( {{\dd q}\over{\dd t}} \big)^2 - V(q) \big] \, \dd t \, $ 
can be defined  with the Simpson quadrature formula  (\ref{simpson})  associated with an
internal polynomial approximation~$ \, q(t) \, $ of degree~2 presented in (\ref{q2t}): 
\moneq \label{LLhh} 
L_h (q_\ell ,\, q_m ,\, q_r) =  {{m\,h}\over{12}} \, \big( g_\ell^2 + 4 \, g_m^2 + g_r^2 \big)
-  {{h}\over{6}} \, \big( V(q_\ell) + 4 \, V(q_m) + V(q_r) \big) . 
\monend 

\smallskip \noindent 
The discrete action $ \, \Sigma_d \, $ for a motion $ \, t \longmapsto q(t) \, $
between the initial time and a given time $ \, T > 0 \, $ is discretized
with $ \, N \, $ regular intervals and take the form 
\moneq \label{SSdd} 
\Sigma_d = \sum_{j=1}^{N-1} L_h (q_j ,\, q_{j+1/2} ,\, q_{j+1} ) .
\monend 
%

%%%%%%%%%%%%%%%%%%%%%%%%%%%%%%%%%%%%%%%%%%%%%%%%%%%%%%%%%%%%%%%%%%%%%%%%%%%%%%%  section 5
\bigskip  \noindent {\bf \large 5) \quad   Discrete Euler-Lagrange equations} 
%%%%%%%%%%%%%%%%%%%%%%%%%%%%%%%%%%%%%%%%%%%%%%%%%%%%%%%%%%%%%%%%%%%%%%%%%%%%%%%%%%%%%%%%%%

\smallskip \noindent
We first write the Maupertuis's stationary-action principle  $ \, \delta \Sigma_d = 0 \, $
with a variation $ \, \delta  q_{j+1/2} \, $ of the internal degree of freedom
in the interval $ \, [t_j ,\, t_{j+1} ] $.
Due to the relations (\ref{derivees-gear})(\ref{derivee-centree}), 
we first observe that
$ \,  {{\partial g_\ell}\over{\partial q_m}} = {{4}\over{h}} $, 
$ \,  {{\partial g_m}\over{\partial q_m}} = 0 \, $ and
$ \,  {{\partial g_r}\over{\partial q_m}} = -{{4}\over{h}} $. 
Then, due to the expression (\ref{LLhh}) of the discrete Lagrangian, 
we have 
$ \,  {{\partial L_h}\over{\partial q_m}} =  {{h}\over{12}} \, \big[ {{8\, m}\over{h}} \,  g_\ell
 - {{8\, m}\over{h}} \,  g_r - 8 \, {{\dd V}\over{\dd q}}(q_m) \big] $.
This partial derivative is equal to zero
when $ \, \delta \Sigma_d = 0 \, $ and
$ \,  m\, {{g_r - g_\ell}\over{h}}  + {{\dd V}\over{\dd q}}(q_m) \, $ 
 is also zero. 
We observe that $ \, g_r - g_\ell = {{4}\over{h}}  \big( q_\ell -2 \,  q_m + q_r \big) \, $
and the condition $ \,  {{\partial L_h}\over{\partial q_m}} = 0 \, $ is finally written
\moneq \label{equation-milieu-intervalle}
m\, {{4}\over{h^2}} \, \big( q_\ell -2 \,  q_m + q_r \big)  + {{\dd V}\over{\dd q}}(q_m) = 0 .
\monend 
We have put in evidence  a second order discretization of the continuous Euler-Lagrange equation
(\ref{oscillateur-nonlineaire})  of this problem. %%%  $ \, {{\dd^2 q}\over{\dd t^2}} +  {{\dd V}\over{\dd q}} = 0 $.
When the hypothesis (\ref{oscillateur harmonique}) 
of an harmonic oscillator is satisfied, we can easily solve this equation and explicit
the middle value $ \, q_m \, $ as a function of the extremities:
\moneq \label{qm-oscillateur-harmonique}
q_m = {{1}\over{1-{{\omega^2 \, h^2}\over{8}}}} \, {{q_\ell + q_r}\over2} .
\monend 
This interpolation is not linear if $ \, h > 0 $. This property illustrates the
underlying polynomial interpolation of degree two. Moreover, 
a stability condition is naturally emerging:
\moneq \label{stabilite}
0 < \omega \, h < 2 \, \sqrt{2} . 
\monend 

\smallskip \noindent
We now incorporate the relation (\ref{qm-oscillateur-harmonique})
inside the expression (\ref{LLhh})  of the discrete Lagrangian.
After a successful formal calculation with the help of the free software
``SageMath'' \cite{sagemath}, we obtain a reduced Lagrangian
\moneq \label{Lagrangien-reduit} 
L_h^r (q_\ell ,\, q_r) = {{1}\over{1-{{\omega^2 \, h^2}\over{8}}}} \Big[
  {1\over2} \, m \, h \, \Big( {{q_r - q_\ell}\over{h}} \Big)^2
  -{{h}\over{2}} \, m\, \omega^2 \, \Big( {{22-h^2 \, \omega^2}\over{48}} \, \big(q_\ell^2 + q_r^2 \big) +
       {1\over12} \, q_\ell \, q_r  \Big) \Big] .
\monend 
The discrete Euler-Lagrange (\ref{Euler-Lagrange-discret}) can now be written
for this reduced Lagrangien (\ref{Lagrangien-reduit}). Instead of the relations (\ref{schema-ordre-deux}), 
we obtain now the following numerical scheme:
\moneq \label{schema-ordre-3} %% \left\{  \begin{array} {c}
{1\over{h^2}} \, \big( q_{j+1} - 2\, q_j + q_{j-1} \big) + {{\omega^2}\over24} \, \big( q_{j+1} +  22 \, q_j + q_{j-1}\big)
- {{\omega^4 \, h^2}\over24} \, q_j  = 0 . 
%% \end{array} \right.
\monend 
The scheme (\ref{schema-ordre-3}) is consistent with the ordinary differential equation (\ref{oscillateur-nonlineaire})(\ref{oscillateur harmonique})
$ \, {{\dd^2 q}\over{\dd t^2}} + \omega^2 \, q(t) = 0 $. 
Secondly, following the definition recalled in  \cite{RM67}, 
the order of truncation   of the scheme (\ref{schema-ordre-3})
is obtained by replacing the discrete variables $\,  q_{j+1} $, $ \,  q_j \, $ and~$ \, q_{j-1} \, $ 
by the solution of the differential equation at the precise points $\, t_j +h $, $ \, t_j \, $ and $ \, t_j -h $. Then 
\moneqstar   \left\{  \begin{array} {l}
 q_{j+1} = q_j + h\, {{\dd q}\over{\dd t}} + {{h^2}\over{2}} \,  {{\dd^2 q}\over{\dd t^2}} + {{h^3}\over{6}} \,  {{\dd^3 q}\over{\dd t^3}}
 + {{h^4}\over{24}} \,  {{\dd^4 q}\over{\dd t^4}} + {{h^5}\over{120}} \,  {{\dd^5 q}\over{\dd t^5}}
 + {{h^6}\over{720}} \,  {{\dd^6 q}\over{\dd t^5}}  + {\rm O}(h^7) \\  \, \vspace{-4 mm} \\ 
q_{j-1} = q_j - h\, {{\dd q}\over{\dd t}} + {{h^2}\over{2}} \,  {{\dd^2 q}\over{\dd t^2}} - {{h^3}\over{6}} \,  {{\dd^3 q}\over{\dd t^3}} 
 + {{h^4}\over{24}} \,  {{\dd^4 q}\over{\dd t^4}} - {{h^5}\over{120}} \,  {{\dd^5 q}\over{\dd t^5}}
 + {{h^6}\over{720}} \,  {{\dd^6 q}\over{\dd t^5}}  + {\rm O}(h^7) .
\end{array} \right. \monendstar 
In these conditions, the left hand side of the relation  (\ref{schema-ordre-3}) is no longer equal to zero  and defines the truncation error
$ \, {\cal T}_h(q_j) $.
With the help of SageMath \cite{sagemath}, one obtains without difficulty the relation
\moneqstar 
{\cal T}_h(q_j) = {1\over1440} \, \omega^6 \, h^4 \, q_j +    {\rm O}(h^6) .
\monendstar 
The numerical scheme  (\ref{schema-ordre-3}) is fourth order accurate  in the sense of the truncation error. 

\smallskip \noindent
A fundamental question concerns stability.  With the linear structure of the finite difference equation (\ref{schema-ordre-3}),
we  consider the equation of degree two obtained by taking
$ \, q_{j-1} = 1 $,\br
$ \, q_{j} = r \, $ and $ \, q_{j+1} = r^2 $.
The scheme is stable when the roots of the corresponding equation are of modulus smaller than 1. 
This equation can we written $ \, a \, r^2 + b \, r + c = 0 \, $
%%%%    c2 =  (hh^2*omega^2 + 24)/24
%%%%    c1 =  (-hh^4*omega^4 + 22*hh^2*omega^2 - 48)/24
%%%%    c0 =  (hh^2*omega^2 + 24)/24
%
with $\, a = 1 + {{h^2 \, \omega^2}\over24} \, $
and 
$ \,  b = -{1\over24} \, \big( 48 - 22 \, h^2\, \omega^2 + h^4\, \omega^4  \big) $.
%%% b = -{2\over{h}} +  {7\over6} \, \omega^2 \, h  -  {5\over32} \, \omega^4 \, h^3  +  {1\over192} \, \omega^6 \, h^5   $.
%   
The discriminant $ \, \Delta \equiv b^2 - 4 \, a \, c \, $ can be factorized: 
\moneqstar  
\Delta = {{\omega^2 \, h^2}\over576} \, (\omega^2 \, h^2 - 24 ) \, (\omega^2 \, h^2 - 12 ) \, (\omega^2 \, h^2 - 8 ) .
\monendstar 

\smallskip \noindent
Under the stability condition (\ref{stabilite}), all the factors in the expression of the discriminant  are negative
and $ \, \Delta < 0 $. Then the equation  $ \, a \, r^2 + b \, r + c = 0 \, $ has two conjugate complex roots $ \, r \, $
and $ \, \overline{r} $. Their product $ \, r \, \overline{r} = | r |^2 \, $ is equal to 1 and the scheme  (\ref{schema-ordre-3})
is stable. 

%%%%%%%%%%%%%%%%%%%%%%%%%%%%%%%%%%%%%%%%%%%%%%%%%%%%%%%%%%%%%%%%%%%%%%%%%%%%%%%  section 6
\bigskip  \newpage \noindent {\bf \large 6) \quad   Symplectic structure} 
%%%%%%%%%%%%%%%%%%%%%%%%%%%%%%%%%%%%%%%%%%%%%%%%%%%%%%%%%%%%%%%%%%%%%%%%%%%%%%%%%%%%%%%%%%

\smallskip \noindent
From the reduced Lagrangian  (\ref{Lagrangien-reduit}), we define  the momentum $ \, p_r \, $ with the
analogue of the relation (\ref{moment}). It comes
\moneq \label{moment-ordre-3} 
%% p_r  = m \, {{q_r - q_\ell}\over{h}} -  m \, {{\omega^2 \, h}\over{12}} \, \big( 7 \, q_r - q_\ell \big) 
%% +   m \, {{\omega^4\,h^3} \over192} \, (15 \, q_r + q_\ell ) -   m \, {{\omega^6\,h^5} \over384} \, q_r \, .
p_r  =  m \, {{q_r - q_\ell}\over{h}} - h \, {{m \, \omega^2}\over6} \,  {{q_\ell + 2 \, q_r}\over{1 - {{\omega^2 \, h^2}\over{8}}}} 
+ h^3 \, {{m \, \omega^4}\over48} \,  {{q_r}\over{1 - {{\omega^2 \, h^2}\over{8}}}} .
\monend

\smallskip \noindent 
This relation (\ref{moment-ordre-3}) can be explicited in the context of  grid points. We have
\moneqstar   \left\{  \begin{array} {l}
  \displaystyle  p_{j+1}  =
  m \, {{q_{j+1} - q_j}\over{h}} - h \, {{m \, \omega^2}\over6} \,  {{q_j + 2 \, q_{j+1}}\over{1 - {{\omega^2 \, h^2}\over{8}}}} 
+ h^3 \, {{m \, \omega^4}\over48} \,  {{q_{j+1}}\over{1 - {{\omega^2 \, h^2}\over{8}}}} 
\\  \vspace{-5 mm} \\
\displaystyle p_j \,\,\,\,\,   =
m \, {{q_j - q_{j-1}}\over{h}} - h \, {{m \, \omega^2}\over6} \,  {{q_{j-1} + 2 \, q_j}\over{1 - {{\omega^2 \, h^2}\over{8}}}} 
+ h^3 \, {{m \, \omega^4}\over48} \,  {{q_j}\over{1 - {{\omega^2 \, h^2}\over{8}}}} .
%%%% m \, {{q_j - q_{j-1}}\over{h}} -  m \, {{\omega^2 \, h}\over{12}} \, (7 \, q_j - q_{j-1} ) 
%%%%  +   m \, {{\omega^4\,h^3} \over192} \, (15 \, q_j + q_{j-1} ) -  m \,  {{\omega^6\,h^5} \over384} \, q_j .
\end{array} \right. \monendstar 

\smallskip \noindent 
We eliminate the variable $ \, q_{j-1} \, $ from these two relations with the help of the 
difference scheme~(\ref{schema-ordre-3}).
We find a recurrence relation for the state $ \, y_j \equiv (p_j ,\, q_j)^{\rm t} $,  
similar to the equation  (\ref{systeme-symplectique-ordre-2}),
but  the matrix $ \, \Phi \, $ is replaced by a matrix $ \, \Phi_3 \, $ that can be explicited:
%
%%%  deno  =  1/24*hh^2*omega^2 + 1
%%%  Phid
%%%  [(hh^4*omega^4 - 22*hh^2*omega^2 + 48)/48   (-hh^5*mm*omega^6 + 36*hh^3*mm*omega^4 - 288*hh*mm*omega^2)/288]
%%%  [           (-hh^3*omega^2 + 8*hh)/(8*mm)                        (hh^4*omega^4 - 22*hh^2*omega^2 + 48)/48]
%%%  fphi11 =  (1/48) * (hh^4*omega^4 - 22*hh^2*omega^2 + 48)
%%%  fphi12 =  (-1/288) * mm * hh * omega^2 * (hh^2*omega^2 - 24) * (hh^2*omega^2 - 12)
%%%  fphi12 =  - mm * hh * omega^2 * (1-hh^2*omega^2/24) * (1-hh^2*omega^2/12)
%%%  fphi21 =  (-1) * 2^-3 * mm^-1 * hh * (hh^2*omega^2 - 8)
%%%  fphi22 =  (1/48) * (hh^4*omega^4 - 22*hh^2*omega^2 + 48)
%
\smallskip
%\moneq \label{Phi3} 
%{\bf \Phi}_3  =  {{1}\over{1 + {{\omega^2  h^2}\over24}}}  
%\left(\!   \begin{array}{cc}   1 \!-\! {11\over24} \,  \omega^2  h^2 +  {{\omega^4  h^4}\over48}  
%  &   \,\,  -m  \, \omega^2 \, h  \, 
%\!  \big(1 - {{\omega^2   h^2}\over12} \big)  \big(1 - {{\omega^2   h^2}\over24} \big)
%\\   {{1}\over{1 - {{\omega^2   h^2}\over8}}} \, {{h}\over{m}}
%&    1 \!-\! {11\over24} \,  \omega^2  h^2 +  {{\omega^4  h^4}\over48}  \end{array} \right)   .
%\monend
%
\moneq \label{Phi3} 
{\Phi}_3  =  {{1}\over{1 + {{\omega^2  h^2}\over24}}}  
\left(\!   \begin{array}{cc}   
			1 \!-\! {11\over24} \,  \omega^2  h^2 +  {{\omega^4  h^4}\over48}  
  			&   \,\,  -m  \, \omega^2 \, h  \, 
\!  \big(1 - {{\omega^2   h^2}\over12} \big)  \big(1 - {{\omega^2   h^2}\over24} \big)
			\\   
			{{h}\over{m}} \, \big({1 \!-\! {{\omega^2   h^2}\over8}}\big)
			&    1 \!-\! {11\over24} \,  \omega^2  h^2 +  {{\omega^4  h^4}\over48}  \end{array} \right)   .
\monend

\smallskip \noindent
We observe that the ``symplectic Simpson'' numerical scheme
defined by  (\ref{systeme-symplectique-ordre-2})(\ref{Phi3})  is an explicit scheme. 
It is easy with SageMath to verify that $ \, {\rm det} \, \Phi_3  = 1 \, $ and
in consequence the scheme is symplectic. 

%%%%%%%%%%%%%%%%%%%%%%%%%%%%%%%%%%%%%%%%%%%%%%%%%%%%%%%%%%%%%%%%%%%%%%%%%%%%%%%  section 7
\bigskip  \noindent {\bf \large 7) \quad   Conservation of a discrete energy} 
%%%%%%%%%%%%%%%%%%%%%%%%%%%%%%%%%%%%%%%%%%%%%%%%%%%%%%%%%%%%%%%%%%%%%%%%%%%%%%%%%%%%%%%%%%

\smallskip \noindent
To explicit a discrete energy that is conserved is not {\it a priori} obvious.
For the harmonic oscillator, we search a conserved quadratic form of the type
\moneq \label{forme-quadratique} 
Q(p,\, q) = {1\over2} \, \xi \, p^2 + \eta \, p \, q +  {1\over2} \, \zeta \, q^2 .
\monend
If we require that $ \, Q(p_{j+1} ,\, q_{j+1}) = Q(p_{j} ,\, q_{j}) \, $
with the  variables $ \, p_{j+1} $, $ \, q_{j+1} $, $ \, p_{j}  \, $ and $ \, q_{j} \, $ 
satisfying a linear dynamics such as (\ref{systeme-symplectique-ordre-2})
with a matrix
\moneqstar
\Phi = \begin{pmatrix} \alpha & \beta \\ \gamma & \delta   \end{pmatrix} 
\monendstar 
of unit determinant, that is    $ \, \alpha \, \delta  -  \beta \, \gamma = 1 $, then the coefficients
$\, \xi $, $ \, \eta \, $ and $ \, \zeta \, $ of the quadratic form (\ref{forme-quadratique})
must satisfy the following homogeneous linear system 
\moneqstar
\Gamma \,  \begin{pmatrix} \xi \\  \eta \\ \zeta  \end{pmatrix} \equiv
\begin{pmatrix} {{\alpha^2-1}\over2} & \alpha \, \gamma & {{\gamma^2}\over2} \\ 
\alpha \, \beta & \alpha \, \delta + \beta \, \gamma - 1 & \gamma \, \delta \\
{{\beta^2}\over2} & \beta \, \delta & {{\gamma^2-1}\over2}   \end{pmatrix}
\,  \begin{pmatrix} \xi \\  \eta \\ \zeta  \end{pmatrix} = \begin{pmatrix} 0 \\ 0 \\ 0 \end{pmatrix} .
\monendstar 
We have
$ \,\, {\rm det} \, \Gamma = (\alpha \, \delta  -  \beta \, \gamma - 1 ) \,  
(1 - \alpha - \delta + \alpha \, \delta  -  \beta \, \gamma ) \, 
(1 + \alpha + \delta + \alpha \, \delta  -  \beta \, \gamma ) \, $
and this expression vanishes when  $ \, \alpha \, \delta  -  \beta \, \gamma = 1 $.
Moreover, when $\, \alpha = \delta $, we obtain $ \, \eta = 0 \, $ and
a conserved quadratic form $ \, Q_c \, $  can be written 
%%% \moneqstar   \label{forme-quadratique-bonne}
$ \, Q_c(p,\, q) = {1\over2} \, \gamma \, p^2 -  {1\over2} \, \beta \, q^2 \, $ 
%% \monendstar 
%
up to a multiplicative constant. Finally, if we set
\moneq \label{hamiltonien-discret}
%%%  modif Juan 17 juillet 2023
H_d \, (p,\, q) \equiv  {1\over{2 \,m}} \, p^2 +  m\, {{\omega^2}\over{2 \,\big(1 - {{\omega^2 \, h^2}\over8} \big)}} \,
  \Big(1 - {{\omega^2 \, h^2}\over12} \Big) \, \Big(1 - {{\omega^2 \, h^2}\over24} \Big) \, q^2 , 
\monend

\smallskip \smallskip  \noindent 
the  symplectic Simpson scheme
\vskip -.4 cm 
\moneqstar \label{systeme-symplectique-ordre-3}
\begin{pmatrix}p_{j+1} \\ q_{j+1} \end{pmatrix} =  \Phi_3 \,  \begin{pmatrix}p_{j} \\ q_{j} \end{pmatrix} 
\monendstar

\smallskip \smallskip  \noindent 
with $ \, \Phi_3 \, $ explicited at the relation (\ref{Phi3}), 
satisfies the following conservation of energy:
\moneqstar
H_d \, (p_{j+1} ,\, q_{j+1}) =  H_d \, (p_{j} ,\, q_{j}) . 
\monendstar 

\newpage
%%%%%%%%%%%%%%%%%%%%%%%%%%%%%%%%%%%%%%%%%%%%%%%%%%%%%%%%%%%%%%%%%%%%%%%%%%%%%%%  section 8
\bigskip  \noindent {\bf \large 8) \quad   First numerical experiments} 
%%%%%%%%%%%%%%%%%%%%%%%%%%%%%%%%%%%%%%%%%%%%%%%%%%%%%%%%%%%%%%%%%%%%%%%%%%%%%%%%%%%%%%%%%%

%
\smallskip \noindent
We have implemented the Simpson symplectic scheme and
have compared it with the Newmark scheme 
(\ref{systeme-symplectique-ordre-2})(\ref{matrice-symplectique-ordre-2}).
Typical results for $ \, N=15 \, $ meshes and one period are displayed on Figure 1.
We have chosen $\, q(t) = \sin \, \omega \, t \, $ and $ \, p(t) = m \, \omega \, \cos \, \omega \, t $, 
with a period $ \, T = 1 $. 
Quantitative errors with the maximum norm are presented in Table~1 below.
An asymptotic order of convergence can be estimated for the momentum, the state and various energies.

\bigskip
%%%%%%%%%%%%%%%%%%%%%%%%%%%%%%%%%%%%%%%%%%%%%%%%%%%%%%%%%%%%%%%%%%%%%%%%%%%%%%%%%%% figure
\smallskip \smallskip %% \vskip -1.3 cm 
%%%  \centerline {\includegraphics[width=.96\textwidth]{simpson-symplectique-article-12fevrier2023.pdf}} 
\centerline {\includegraphics[width=.96\textwidth]{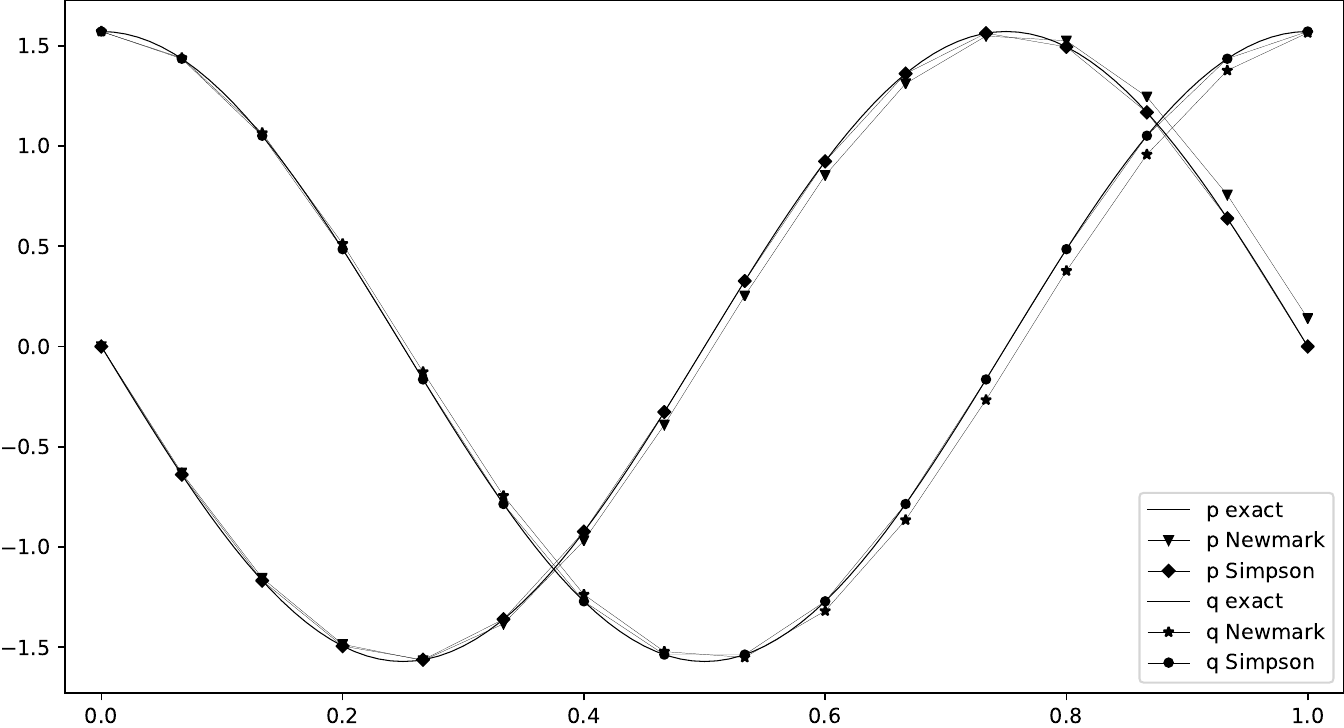}} 

\smallskip  \noindent
Figure 1. Typical evolution of an harmonic oscillator.
The momentum $ \, p \, $ follows a cosine curve and the state $ \, q \, $ a sine shape. 
Comparison of the exact solution and the Newmark and  symplectic Simpson schemes for $ \, N = 15 \, $ meshes.
Both schemes give very satisfactory results. Observe that the momentum $ \, p \, $ is
very close to the exact solution with the  symplectic Simpson scheme.
Observe that the momentum data have been rescaled. 
%%%%%%%%%%%%%%%%%%%%%%%%%%%%%%%%%%%%%%%%%%%%%%%%%%%%%%%%%%%%%%%%%%%%%%%%%%%%%%%%%%%

\bigskip \smallskip \smallskip
%%%%%%%%%%%%%%%%%%%%%%%%%%%%%%%%%%%%%%%%%%%%%%%%%%%%%%%%%%%%%%%%%   table 
{\centering{\begin{tabular}{|c|c|c|c|c|c|c|c|c|c|c|}    \hline
& number of meshes & 10 & 20 & 40 & $\,$ order $\,$ \\   \hline
Newmark & momentum & $ 1.91 \,\, 10^{0} $  &  $ 5.02 \,\, 10^{-1} $  & $ 1.27  \,\, 10^{-1} $  & 2 \\   \hline
Symplectic Simpson & momentum & $  3.41 \,\, 10^{-3} $  &  $  2.17 \,\, 10^{-4} $   &  $  1.35  \,\, 10^{-5} $        & 4 \\   \hline
Newmark & state    & $ 2.38  \,\, 10^{-1} $  &  $ 6.12  \,\, 10^{-2} $  & $ 1.55   \,\, 10^{-2} $  & 2  \\   \hline
Symplectic Simpson  & state    & $ 4.11  \,\, 10^{-4} $  & $ 2.55 \,\, 10^{-5} $   & $ 1.60 \,\, 10^{-6} $     & 4  \\   \hline
Newmark & energy (\ref{HH-exact}) &  $ 1.70 \,\, 10^{-13} $ &   $ 4.97 \,\, 10^{-14} $  & $ 1.49 \,\, 10^{-13} $   & exact  \\   \hline
Symplectic Simpson  & energy (\ref{HH-exact}) &  $ 2.51  \,\, 10^{-2}  $ &  $ 1.67 \,\, 10^{-3}  $ &  $ 1.03  \,\, 10^{-4}  $ & 4  \\   \hline
Newmark & energy (\ref{hamiltonien-discret}) &   $ 2.70  \,\, 10^{-2} $  &   $ 1.67 \,\, 10^{-3} $  &  $ 1.03 \,\, 10^{-4} $   & 4  \\   \hline
Symplectic Simpson  & energy (\ref{hamiltonien-discret}) &  $ 2.13 \,\, 10^{-14}  $  & $ 2.49  \,\, 10^{-13}  $  & $ 1.63 \,\, 10^{-13}  $ & exact
\\   \hline \end{tabular}}}
%%%%%%%%%%%%%%%%%%%%%%%%%%%%%%%%%%%%%%%%%%%%%%%%%%%%%%%%%%%%%%%%%%%%%%%%%%%%%%%%%%%%%%%%%%

\smallskip  \smallskip \smallskip \noindent
Table 1. Errors in the %% $ \, \ell^\infty \, $
maximum norm.
We observe again that the momentum is very well approximated with the symplectic Simpson scheme. 
The estimated order of convergence is the closest integer $ \, \alpha \, $
measuring  the ratio of successive errors in a given line by a negative power of 2
of the type~$ \, 2^{-\alpha} $.

%%%%%%%%%%%%%%%%%%%%%%%%%%%%%%%%%%%%%%%%%%%%%%%%%%%%%%%%%%%%%%%%%%%%%%%%%%%%%%%  section 9
\bigskip  \noindent {\bf \large 9) \quad    Conclusion and perspectives} 
%%%%%%%%%%%%%%%%%%%%%%%%%%%%%%%%%%%%%%%%%%%%%%%%%%%%%%%%%%%%%%%%%%%%%%%%%%%%%%%%%%%%%%%%%%

\smallskip \noindent
The symplectic Simpson numerical scheme has been developed in this contribution.
It has been tested for an harmonic oscillator.
The method is symplectic, conditionnally stable and
is fourth order accurate for state and momentum.

\smallskip \noindent
An important question is still open concerning the nonlinear case.
The elimination of the internal degree of freedom is not possible in that case.
Once this question has a satisfactory answer, the extension to systems with mutiple degrees of freedom 
is a natural objective for future studies.

%%%%%%%%%%%%%%%%%%%%%%%%%%%%%%%%%%%%%%%%%%%%%%%%%%%%%%%%%%%%%%%%%%%%%%%%%%%%%
%%%%%%  \bigskip     \noindent {\bf \large     Acknowledgments} 
%%%%%%%%%%%%%%%%%%%%%%%%%%%%%%%%%%%%%%%%%%%%%%%%%%%%%%%%%%%%%%%%%%%%%%%%%%%%%%%%%%%%%%%%%%

%%%%%%  \noindent
%%%%%%  The author is particularly grateful to Juan Antonio Rojas-Quintero for very constructive discussions preliminary to this work.

%%%%%%%%%%%%%%%%%%%%%%%%%%%%%%%%%%%%%%%%%%%%%%%%%%%%%%%%%%%%%%%%%%%%%%%%%%%%  references
\bigskip       \noindent {\bf  \large  References }
%%%%%%%%%%%%%%%%%%%%%%%%%%%%%%%%%%%%%%%%%%%%%%%%%%%%%%%%%%%%%%%%%%%%%%%%%%%%%%%%%%%%%%%%%%

 %%%%%%%%  \vspace{-.24cm}

\end{document}